\documentclass[12pt]{article}
\usepackage{amsmath}
\usepackage{amssymb}
\usepackage{amscd}
\usepackage{amsfonts}

\def\reel{\hbox{{\rm R}\kern-1em\hbox{{\rm I} }}}
\def\relatif{\ \hbox{{\rm Z}\kern-.4em\hbox{\rm Z}}}
\def\nat{\hbox{{\rm N}\kern-1em\hbox{{\rm I} } }}
\def\comp{\hbox{{\rm C}\kern-.55em\hbox{{\rm I} } }}
\def\smallcomp{\hbox{\fiverm C}\kern-.35em{\hbox{\fiverm I}}}

\def\fudge{\mathchoice{}{}{\mkern.5mu}{\mkern.8mu}}
\def\bbc#1#2{{\rm \mkern#2mu\vbar\mkern-#2mu#1}}
\def\bbb#1{{\rm I\mkern-3.5mu #1}} \def\bba#1#2{{\rm #1\mkern-#2mu\fudge
#1}}
\def\bb#1{{\count4=`#1 \advance\count4by-64 \ifcase\count4\or\bba
A{11.5}\or \bbb B\or\bbc C{5}\or\bbb D\or\bbb E\or\bbb F \or\bbc G{5}\or\bbb H\or \bbb I\or\bbc J{3}\or\bbb
K\or\bbb L \or\bbb M\or\bbb N\or\bbc O{5} \or \bbb P\or\bbc Q{5}\orrrr\b bb R\or\bbc S{4.2}\or\bba T{10.5}\or\bbc
U{5}\or    \bba V{12}\or\bba W{16.5}\or\bba X{11}\or\bba Y{11.7}\or\bba Z{7.5}\fi}}
\def\rat{\hbox{{\rm Q}\kern-.70em\hbox{{\rm I} } }}

\newcommand{\refm [1]} {(\ref{#1})}

\newcounter{corcountrer}
\newcounter{theoremcounter}
\newcounter{lemmacounter}

\newcounter{remarkcounter}
\newcounter{propositioncounter}

\newtheorem{cor}[corcountrer]{Corollary}
\newtheorem{thm}[theoremcounter]{Theorem}
\newtheorem{lemma}[lemmacounter]{Lemma}

\newtheorem{rem}[remarkcounter]{Remark}
\newtheorem{prop}[propositioncounter]{Proposition}

\newcommand{\be}{\begin{equation}}
\newcommand{\ber}{\begin{eqnarray}}

\newcommand{\nin}{\noindent}
\newcommand{\non}{\nonumber}

\newcommand{\un}{\underline}
\def\qed{\hfill \vrule height1.3ex width1.2ex depth-0.1ex}
\def\bbb#1{{\rm I\mkern-3.5mu #1}} \def\bba#1#2{{\rm #1\mkern-#2mu\fudge
#1}}

\newcommand{\la}{\label}

\title{Clustering in
coagulation - fragmentation processes, random combinatorial structures and additive number systems: Asymptotic
formulae and limiting laws.}

\author{{\bf Gregory A. Freiman}
\thanks{E-mail: grisha@math.tau.ac.il }\\
School of Mathematical Sciences, Raymond and Beverly Sackler Faculty\\
of Exact Sciences, Tel-Aviv University, Ramat-Aviv,\\ Tel-Aviv,
Israel.\\
\quad {\bf Boris L. Granovsky}
\thanks{E-mail: mar18aa@techunix.technion.ac.il} \\
Department of Mathematics, Technion-Israel Institute of Technology,\\
Haifa, 32000, Israel.}
\begin{document}
\maketitle \vskip 5cm

\nin American Mathematical Society 2000 subject classifications:

\nin Primary-60K35, 05A15; secondary-05A16, 05C80, 11M45.

\nin Keywords and phrases: Coagulation-fragmentation process,  Random combinatorial structures, Local limit theorem,
Distributions on the set of
partitions, Additive number systems.

\newpage

\begin{center} {\bf Abstract}
\end{center}

 \nin We develop a unified approach to the problem of clustering
in the three different fields of applications indicated in the
title of the paper, in the case when the parametric function of
the models is regularly varying with positive exponent.  The
approach is based on Khintchine's probabilistic method that grew
out of the Darwin-Fowler method in statistical physics. Our main
result is the derivation of asymptotic  formulae for the
distribution of the largest and the smallest clusters (=
components), as the total size of a structure (= number of
particles) goes to infinity. We discover that  $n^{\frac{1}{l+1}}$
is the threshold for the limiting distribution of the largest
cluster. As a by-product of our study, we
 prove   the
 independence of the numbers of groups of fixed sizes,  as $n\to \infty.$
  This is  in
 accordance with the general principle of  asymptotic independence of
 sites in mean-field models. The latter principle is commonly accepted
 in statistical physics, but
 not rigorously proved.
\newpage

\nin \section{Introduction: The objective and the context.} \setcounter{equation}{0}

\nin We develop a unified approach to the problem of clustering in the three different fields of applications
indicated in the title of the paper. The approach is based on  Khintchine's probabilistic method that grew out of
the Darwin-Fowler method in statistical physics. To the best of our knowledge, the first application of
Khintchine's method  for coagulation - fragmentation processes was made in
 \cite{frgr}, where it was used for the derivation of asymptotic formulae
 for the partition function of the invariant measure of the process.
The present paper extends the method to  much more complicated asymptotic problems arising in the study of
clustering.

\nin  Our main
result is the derivation of asymptotic  formulae for the distribution of the largest and the smallest clusters
(= components), as the total size of a structure (= number of particles) goes
to infinity.

\nin The organization of the paper is as follows. Section 2 provides a formal  mathematical setting that
encompasses the clustering problems arising in the contexts of coagulation - fragmentation processes, random
combinatorial structures and additive number systems. The mathematical problem is stated as follows. Let the
functions $g,S:\comp\rightarrow \comp $ be related via $g(z)=e^{S(z)}, \ \vert z\vert <R, R>0$. Under a given
asymptotic behavior of the Taylor coefficients $\{ a_n\}$ of the function $S$ one must  explore the asymptotic
behavior of certain quantities related to
 the Taylor coefficients
 $\{c_n\}$ of the function $g.$

\nin  The problem is considered for the class
 of functions $S,$ such that
$a_n\sim n^{l-1}L(n), \ l>0, \ n\to \infty,$ where $L$ is a slowly varying function. A specific feature of this
class of functions is that it provides the validity of the normal local limit theorem for the associated
probabilistic model.

\nin In Section 3 we explain  the idea of Khintchine's method and apply it to the  derivation of  the   asymptotic
formulae for  the limiting distributions of the
largest and the smallest clusters. We find that $n^{\frac{1}{l+1}}$ is the threshold for the limiting distribution
of the largest cluster.

\nin In Sections 4-6 we demonstrate how to interpret these asymptotic formulae in the
 context of the aforementioned three  fields, and  we provide a description of
the striking picture of clustering that follows from these formulae.
It turns out that for large $n$, almost all weight of $n$
 is distributed  into groups of sizes about $n^{\frac{1}{l+1}}$, while the rest of the weight is
 made up of  groups of small sizes.
\nin As a by-product of our study:

\nin (i) We
 prove   the
 independence of the numbers of groups of fixed sizes,  as $n\to \infty.$
  This is  in
 accordance with the general principle of  asymptotic independence of
 sites in mean-field models. The latter principle is commonly accepted
 in statistical physics, but
 not rigorously proved.

\nin (ii) We recover an  asymptotic result by J. Knopfmacher, A. Knopfmacher and
 R. Warlimont, that is widely
known in the theory of additive number
systems.

\nin \section{Mathematical setting and preliminaries.} \setcounter{equation}{0} \nin We consider throughout the
paper the set  ${\cal F}(l),  l> 0,$ of
  sequences $a=\{a_n\}_1^\infty, \quad a_n\ge 0,\  n\ge 1,$
 with the following
asymptotic behavior:

\be a_n\sim n^{l-1}L(n),\quad \mbox{as $n\to \infty$},
 \quad l> 0. \la{x7} \end{equation}

\nin Here and in what follows,
 $L$ is  a slowly varying (s.v.) function  at infinity
 (for references see
 \cite{bing},
   \cite{sen}).
\nin We will need the  following two asymptotic properties of s.v. functions:

\be L(x)=o(x^\epsilon),\quad  \mbox{as $x\to \infty$}, \quad \mbox{for all}\quad  \epsilon>0, \la{L1}
\end{equation}

\be x^{-\epsilon}=o(L(x)),\quad  \mbox{as $x\to \infty$}, \quad \mbox{for all}\quad  \epsilon>0.
\la{L2}\end{equation}

\nin We assume further that
\begin{itemize}

\item $L$ is  differentiable  on $[0,\infty)$. This is
based on the
 fact (\cite{sen}, p.17)
 that for any
s.v. function $L$ there exists a s.v. function
$\tilde{L},$  that possesses the aforementioned property and satisfies
$L(x)\sim \tilde{L}(x), $ as $x\to \infty.$

\item The function
$x^{-\delta}L(x)$ is locally bounded on $[0,\infty)$, for any $\delta>0.$
\end{itemize}

\nin It is easy  to derive from the   representation of the set of s.v. functions (\cite{sen}, p.2) that  the
sequences $a\in {\cal F}(l), \ l>0,$ satisfy

\be
 \lim_{n \to \infty} \frac{a_n}{a_{n+1}}=1.
\label{rs}
\end{equation}

\nin We will also need the fact  that a s.v. function $L$ has  a conjugate function $L^*$  (\cite{sen}, p. 25 and
\cite{bing}, p.47), which is also  a s.v. function and is uniquely defined (up to asymptotic equivalence)  by the
asymptotic relationship \be L^*(x)L(xL^*(x))\sim L(x)L^*(xL(x))\sim 1,\quad {as} \quad x\to \infty. \la {cnj}
\end{equation}

\nin  (\ref{cnj}) says that  the asymptotic behavior of $L^*$ is converse to the one of $L,$  in the sense that
\be \lim_{n\to \infty} L^*(n)=\big(\lim_{n\to \infty} L(n)\big)^{-1}, \la{cn1}\end{equation}

\nin provided the limits exist. \nin The characterization of the class of sequences ${\cal F}(l), \ l>0,$ is given
by the celebrated Karamata Tauberian theorem (for references see \cite{sen}, p.59 and \cite{fel}, p.423) that is a
widely used tool in different fields of probability.

\nin In effect, we will employ the following corollary of Karamata's theorem.

\vskip.5cm
\begin{thm}\label{vf}(\cite{fel}, p.423)

\nin Let a sequence $a=\{a_n\ge 0, \ n\ge 1\}$ be  ultimately
monotone, and suppose that  the radius of convergence of the power
series (in $z$) \be S(z) = \sum_{n=1}^{\infty} {a_n z^n},\quad
z\in \comp \label{gen}\end{equation} \nin equals 1.
 Then the  two conditions (i) and (ii)
 are
equivalent:

\nin (i)

\be S(z)\sim \frac{\Gamma(l)}{(1-z)^l}L\Big(\frac{1}{1-z}\Big), \quad l>0, \quad \mbox{as $ z\to 1^{-}$}, \la{kar}
\end{equation} \nin where $\Gamma$  is the gamma function, and

 \nin  (ii) \be a_n \sim n^{l-1}L(n)\in{\cal F}(l), \quad l> 0.
 \end{equation}

\end{thm}

\vskip 1cm

\nin For our subsequent study we will make use only of the abelian part (i) of the above theorem.

\nin Next we define the sequence $c=\{c_n\}_0^\infty $ generated by the above sequence $a$ in the following
manner: \be g(z):=\sum_{n=0}^{\infty} c_n z^n=e^{S(z)}, \quad \vert z\vert<1. \la{00}\end{equation}

\nin We will demonstrate in Sections 4-6 that the above  form of the exponential relationship between two
generating functions arises in the three fields in the title of the present  paper. In view of this,  a variety of
problems related to \refm[00] (but quite different from the problem considered by us)  have been studied by many
researchers.

\nin  Based on (\ref{rs}),   it is easy to derive (see \cite{dgg} and \cite{bur},Lemma 1.22) that the radius of
convergence of  the series for $g(z)$ equals $1.$
 Moreover, it was recently proven by J. Bell and S. Burris
(\cite{bbur}, Lemma 4.2) that (\ref{rs}) implies

\be
 \lim_{n \to \infty} \frac{c_n}{c_{n+1}}=1.
\label{rs1}
\end{equation}

\nin This fact is important, since by  Compton's density theorem
(see for references \cite{bbur} and  \cite{bur}, Ch.4), the
condition \refm[rs1] implies that all partition sets of an
additive number system  have asymptotic density which is either
$0$ or $1$.

\nin To formulate the problem of clustering that is  addressed in the present paper we introduce some more
notation.
 For  given $r, n,\  1\le r= r(n)\le n, \ n=1,2,\ldots,$
 we denote
\be \underline{S}_n^{(r)}(z)=
 \sum_{j=1}^{r} {a_j z^j},  \qquad
\bar{S}_n^{(r)}(z)=\sum_{j=r}^{n} {a_j z^j},  \quad \vert z\vert< 1 \la{gen1}
 \end{equation}
\nin and  consider the two power series

\be \underline{g}_n^{(r)}(z)=e^{\underline{S}_n^{(r)}(z)}: =\sum_{j=0}^{\infty}\underline{ c}_j^{(r)} z^j \quad
{and}\quad \bar{g}_n^{(r)}(z)=e^{\bar{S}_n^{(r)}(z)}:=\sum_{j=0}^{\infty} \bar{c}^{(r)}_j z^j, \quad \vert z\vert<
1. \la{gen2}
\end{equation}


\nin  Setting $r=n^\beta, \
 0\le \beta<1$ and denoting $\underline{c}_n^{(n)}=c_n, \ n=1,2,\ldots,$
 our ultimate objective will be the derivation of the
limits, as $n\to \infty,$ for the two quantities \be
\underline{d}^{(r)}_n:=\frac{\underline{c}^{(r)}_n}{c_n}, \quad
{and}\quad \bar{d}^{(r)}_n:=\frac{\bar{c}^{(r)}_n}{c_n}. \la{x121}
\end{equation}

\nin Here and in what follows we agree that $r=\bullet$ means that
$r=[\bullet],$ where $[\bullet]$ is the integer part of the
 number $\bullet$.

 \nin  \section{Asymptotic formulae and limiting laws}
\setcounter{equation}{0}

 \nin   We will study the  above posed
problem with the help of the probabilistic method formulated by
Khintchine in \cite{Kh}, Ch.IV,V (see also \cite{frgr}).
Independently of the context of the problem considered, the
implementation of Khintchine's method
 for deriving  asymptotic formulae always follows the following
 two - step scheme:

\nin(i) The construction of
 an auxiliary probabilistic
 model with a free parameter
 that enables one  to express the quantity in question via the probability
 function of a
 sum of independent integer-valued random variables forming a triangular
 array.

\nin (ii) The proof of the  normal local  limit theorem via a proper choice
 of a free parameter in the probabilistic model in (i).

\nin The problem formulated in (\ref{x121}) requires the
derivation of asymptotic formulae  for the coefficients
$\underline{c}^{(r)}_n$ and $\bar{c}^{(r)}_n,$ for all $r=n^\beta,
\ 0\le\beta< 1.$ In the case $L(x)\equiv 1$ such a formula for
$\underline{ c}^{(n)}_n$ was established in \cite{frgr}, with the
help of Khintchine's method.
 Our primary aim in this section will
be to extend the method  to the aforementioned problem (\ref{x121}). This will
require a much more complicated asymptotic analysis.

\nin The  probabilistic model suggested below  is a modification of the one
 in \cite{frgr}.
 We start by setting in (\ref{gen1}), (\ref{gen2})
\be z=e^{-\sigma +2\pi i\alpha}, \label{-1}
\end{equation}
for some  $\sigma,\alpha\in R.$ Then, analogous to Lemma 1 in \cite{frgr},
  the following representations
of $\underline{ c}^{(r)}_n,$ and $ \bar{c}^{(r)}_n$ are valid: $$ \underline{c}_n^{(r)}= e^{n\sigma}\int _0^1
\prod_{j=1}^{r}\left( \sum_{k=0}^{\infty}\frac{a_j^k e^{-jk\sigma+2\pi i\alpha jk}}{k!}\right)\times e^{-2\pi
i\alpha n}d\alpha, $$ $$ \bar{c}_n^{(r)}= e^{n\sigma}\int _0^1  \prod_{j=r}^{n}\left(
\sum_{k=0}^{\infty}\frac{a_j^k e^{-jk\sigma+2\pi i\alpha jk}}{k!}\right)\times e^{-2\pi i\alpha n}d\alpha, $$ \be
1\le r\le n, \quad n=1,2,\ldots, \la{cn}\end{equation} \nin where $\sigma\in R$ is arbitrary. For this reason
$\sigma$ is called a free parameter.
 It plays an important role in  the method.

\nin To attribute a probabilistic meaning to the RHS's in
(\ref{cn}), we make use of the following notation:

\be p_{jk}=\frac{(a_{j} e^{-\sigma j})^k}{k!\exp\left( a_je^{-\sigma j}\right)}, \quad j=1,\ldots,n,\quad
k=0,1,\ldots \label{p} \end{equation}

\be \varphi_j(\alpha)=  \sum_{k=0}^{\infty} p_{jk}e^{2\pi i\alpha jk}, \quad \alpha\in R , \quad 1\le j\le n,
\label{v1} \end{equation}

\be \underline{\varphi}^{(r)}(\alpha)= \prod_{j=1}^r \varphi_j(\alpha) , \quad \bar{\varphi}^{(r)}(\alpha)=
\prod_{j=r}^n \varphi_j(\alpha), \quad \alpha\in R . \label{v2} \end{equation}

\nin Notice  that for a given $j,$ (\ref{p}) can be viewed as the Poisson probability function with parameter
$a_je^{-\sigma j},\ \sigma\in R.$
 Now (\ref{cn})  can be rewritten as \be \underline{c}_n^{(r)}=
\exp\Big(\underline{S}_n^{(r)}(e^{-\sigma})+n\sigma\Big)\int_0^1 \underline{\varphi}^{(r)}(\alpha)e^{-2\pi i\alpha
n}d\alpha, \la{3y}
\end{equation}

\be \bar{c}_n^{(r)}= \exp\Big(\bar{S}_n^{(r)}(e^{-\sigma})+n\sigma\Big)\int_0^1
\bar{\varphi}^{(r)}(\alpha)e^{-2\pi i\alpha n}d\alpha,\quad 1\le r\le n, \quad n=1,2,\ldots \label{3y1}
\end{equation}

\nin  The  representations (\ref{3y}), (\ref{3y1})  belong to  the
core of Khintchine's method. The idea behind the  representations
is that $\un{\varphi}^{(r)}(\alpha)$ in (\ref{3y}) can be
interpreted as a characteristic function of the sum
$\underline{Y}_{n}^{(r)}=X_1+\ldots+X_r $ of independent lattice
random variables $X_1,\ldots, X_r, \  1\le r\le n$,
 defined by

\be Pr(X_j=jk)=p_{jk}, \quad j=1,\ldots, r , \quad k=0,1,\ldots \label{pr}
\end{equation}

\nin Hence,

\be \underline{T}_n^{(r)}=\underline{T}_n^{(r)}(\sigma):=\int_{0}^{1} \underline{\varphi}^{(r)}(\alpha)e^{-2\pi
i\alpha n}d\alpha = Pr(\underline{Y}_{n}^{(r)}=n). \label{333} \end{equation}

\nin Analogously, writing $\bar{Y}_{n}^{(r)}=X_r+\ldots+X_n, \ 1\le r\le n,$ we get

\be \bar{T}_n^{(r)}=\bar{T}_n^{(r)}(\sigma):= \int_{0}^{1}\bar {\varphi}^{(r)}(\alpha)e^{-2\pi i\alpha n}d\alpha =
Pr(\bar{Y}_{n}^{(r)}=n). \label{3333} \end{equation}

\nin In view of (\ref{3y}), (\ref{3y1}) and (\ref{333}), (\ref{3333}), we will focus now on finding the asymptotic
behavior of the probabilities $P(\underline{Y}_{n}^{(r)}=n)$ and $P(\bar{Y}_{n}^{(r)}=n),$ as $n\to \infty.$

\nin First recall that  the classical normal local limit theorems
(see \cite{gn}, \cite{po}, p.78, \cite{il}) are restricted to the
case of a sum of  independent  lattice random variables, while in
our case, as we will see later on,  the lattice random variables
$X_j$ given by \refm[pr] with $\sigma$ depending on $n$ form a
triangular array.
 So, even the existence of the limiting probability
density for our problem is in question.

\nin Notwithstanding  this, we will demonstrate (Theorem 1 below)
that   a proper choice of the free parameter $\sigma$
 guarantees a version of  the famous Gnedenko local limit theorem.

\nin Let  $\underline{M}^{(r)}_n=\un{M}^{(r)}_n(\sigma):= E\underline{Y}^{(r)}_n,$\ \ $(\underline
B^{(r)}_n)^2=(\underline B^{(r)}_n)^2(\sigma):=Var \underline{Y}^{(r)}_n $ and
$\underline\rho^{(r)}_n=\underline\rho^{(r)}_n(\sigma): =E(Y^{(r)}_n-EY^{(r)}_n)^3,$ and denote by
$\bar{M}^{(r)}_n=\bar{M}^{(r)}_n(\sigma)$, $(\bar{B}^{(r)}_n)^2=(\bar{B}^{(r)}_n)^2(\sigma)$ and
$\bar{\rho}^{(r)}_n=\bar{\rho}^{(r)}_n(\sigma)$ the same moments of the sum $\bar{Y}^{(r)}_n.$
 \nin It
follows from (\ref{pr}) that $j^{-1}X_j, \ j=1,\ldots, n,$ are Poisson($a_je^{-\sigma j}$) random variables. So,
we have the following expressions for the above quantities:
 \be
\underline{M}^{(r)}_n=\sum_{j=1}^r ja_j e^{-\sigma j},\quad \bar{M}^{(r)}_n=\sum_{j=r}^n ja_j e^{-\sigma j}, \quad
n=1,2,\ldots \la{gne0} \end{equation} \be
 (\underline B^{(r)}_n)^2
  = \sum_{j=1}^r j^2a_j e^{-\sigma j}, \quad (\bar{B}^{(r)}_n)^2
  =\sum_{j=r}^n j^2a_j
 e^{-\sigma j},
\quad n=1,2,\ldots. \la{gne1}
\end{equation}
\be \underline\rho^{(r)}_n =\sum_{j=1}^r j^3a_j e^{-\sigma j},\quad \bar{\rho}^{(r)}_n=\sum_{j=r}^n j^3a_j
e^{-\sigma j}, \quad n=1,2,\ldots. \la{gne2}
\end{equation}

\nin Now we choose in (\ref{3y}) (resp.(\ref{3y1})) the parameter $\sigma $ equal to the unique solution of the
equations (\ref{sig}) (resp.(\ref{sig1}) below:

 \be \underline{M}^{(r)}_n(\sigma)=n \la{sig} \end{equation}
\nin and

 \be \bar{ M}^{(r)}_n(\sigma)=n. \la{sig1} \end{equation}

\nin The existence and  uniqueness of the solution (in $\sigma$)
of each of the two equations (\ref{sig}), (\ref{sig1}), for any
given $1\le r\le n$ and  $n=1,2,\ldots, $ follow from the
assumption that $a_j>0,\ j=1,2,\ldots.$

\nin The idea of the  above   choice of the free parameter $\sigma $, that goes back to Khintchine's book
\cite{Kh}, is to evaluate the  probabilities  in \refm[333],
 \refm[3333] when $n$ is the ``most probable value"  of
the sums $\underline{Y}_{n}^{(r)}, $ $\bar{Y}_{n}^{(r)}$. This makes the exponential factor in the expression of
the  normal density equal to $1$, which  will enable us to obtain the principal term in the asymptotic expansions
for the above probabilities .
  We will assume further on  that $a\in{\cal F}(l),$  $ l>0$
and denote by $\underline\sigma_n^{(r)},$ $\bar{\sigma}_n^{(r)}$ the solutions of (\ref{sig}), (\ref{sig1})
correspondingly.

\nin \begin{rem}\label{89} In statistical physics, the  idea of
introducing a free parameter has its roots in the famous
Darwin-Fowler asymptotic method developed in the 1920s for
evaluating partition functions and mean values of occupation
numbers. A good exposition of the method is given in \cite{ehse},
Ch.6. In this method, the above quantities are expressed as
complex integrals over a circle around the origin, of an arbitrary
radius (= free parameter). Evaluating the integrals by the method
of steepest descents, the free parameter is taken to be equal to
the unique minimum
  point   in $[0,1]$ of the  integrand.
In the preface to  his book \cite{Kh} Khintchine writes that the
main novelty of his approach consists of replacing ``the
complicated analytical apparatus (the method of Darwin-Fowler) by
.....the well developed limit theorems  of the theory of
probability.....that can form the analytical basis for all the
computational formulas of statistical physics."

\nin Finally, notice that a probabilistic method for the study of
asymptotic problems arising in enumeration of permutations was
quite independently suggested in the 1940s by V. Goncharov.
Subsequently the method  was extensively developed by generations
of researchers who applied it to general combinatorial structures.
The history of this line of research can be found in Kolchin's
book \cite{kol}.
\end{rem}
\begin{rem}\label{67} As we already mentioned, a specific feature
 inherent in
Khintchine's method is that  the free parameter $\sigma$ depends
on $n$, so that the random variables  $X_j,\ j=1,2,\ldots$ form a
triangular array. In the case of an array the conditions for a
normal local limit theorem are not known. For this reason,
starting from A. Khintchine (see \cite{Kh}, Ch.IV) and until the
present time, the establishment of a local limit theorem for sums
of random variables  depending on a free parameter required
sophisticated asymptotic analysis that differed from problem to
problem. As examples, see (in  chronological order)\cite{fr} of G.
Freiman, \cite{po},Ch.2 of A. Postnikov, \cite{frp} of G. Freiman
and J. Pitman, \cite{mu} of R. Mutafchiev, \cite{kol} of V.
Kolchin, \cite{frgr} of G. Freiman and B. Granovsky and \cite{fvy}
of G. Freiman, A. Vershik and Yu. Yakubovitz. In particular, note
that \cite{dfm} of J. Deshouillers, G. Freiman and W. Moran gives
an example of an array of random variables  for which the local
limit theorem fails though the Lyapunov condition holds.
\end{rem}

\nin Throughout the paper we will denote by $h,\ h_i, \ i=1,2,
\ldots$ positive constants that appear in asymptotic formulae.

\nin The following basic property of the solutions $\un{\sigma}_n,
\ \bar{\sigma}_n$ allows for the implementation of \refm[kar]. The
proof of it is similar to that of Lemma 3 in \cite{frgr}.

   \nin \begin{lemma}\label{xz}

\nin Let $n\ge r\ge n^\epsilon$, for some $\epsilon>0.$ Then

\be \lim _{n\to \infty}\underline{\sigma}_n^{(r)}= 0, \quad \lim _{n\to \infty}\bar{\sigma}_n^{(r)}= 0.
\la{4}\end{equation}
\end{lemma}

\nin It is clear that the straightforward application of the
summation formula \refm[kar] to the sums in
\refm[gne0]-\refm[gne2] is not possible. Our subsequent asymptotic
analysis extends the one in  \cite{frgr} in two different
directions:   from $c_n$ to
$\underline{c}_n^{(r)},\bar{c}_n^{(r)}, \ r=n^\beta, \ 0\le
\beta\le 1, $ and from the smooth case $a_n\sim n^{l-1}, \ l>0$ to
the case $a_n\sim n^{l-1}L(n), \ l>0.$ Our main tools will be the
Euler integral test and a summation theorem of Abelian type.

\nin  Consider the function $f(x,\sigma)=x^l L(x)e^{-\sigma x},\ x>0, \ \sigma\in R,\ l>0.$ \nin If $\sigma>0,$
then
 for sufficiently large $x>0$ and sufficiently small $\sigma$
 the function
 $f$ has a  maxima at the point $x=x(\sigma)$ which is
 the solution of the equation
\be xL'(x)+ (l-\sigma x)L(x)=0. \la{Ls}
\end{equation}
\nin Since (see \cite{sen},p.7) \be \lim_{x\to \infty}\frac{xL'(x)}{L(x)}=0, \la{mn}
\end{equation}
\nin for any  s.v. function $L,$ the asymptotic solution of
(\ref{Ls}) is given by $x \sim l\sigma^{-1},$ as  $\sigma\to 0^+$.
 In the case $\sigma\le 0,$ the function $f$ is increasing in $x$
for sufficiently large $x>0.$
 Since we are interested in $r=n^\beta,\ 0<\beta<1, $ Lemma 1 is
valid. So,  applying in  both cases of $\sigma$ the integral test
to the sums $\underline{M}_n{(r)}, \ \bar{M}_n{(r)}$, we can
rewrite (\ref{sig}) and (\ref{sig1}) as

$$ n= \underline{M}_n{(r)}\sim \int_{ 1}^r
f(x,\underline{\sigma}_n^{(r)})dx= $$
\be\Big(\vert\underline{\sigma}_n^{(r)}\vert\Big)^{-(l+1)}
\int_{\vert\underline{\sigma}_n^{(r)}\vert}^{r\vert\underline{\sigma}_n^{(r)}\vert}
t^lL\Big(\frac{t}{\vert\underline{\sigma}_n^{(r)}\vert}\Big)
\exp\Big(-tsign(\underline{\sigma}_n^{(r)})\Big)dt,\quad
 l>0, \quad n\to \infty \la{21}\end{equation}
 \nin and

$$ n= \bar{M}_n{(r)}\sim \int_{ r}^n
f(x,\bar{\sigma}_n^{(r)})dx=$$
\be
\Big(\vert\bar{\sigma}_n^{(r)}\vert\Big)^{-(l+1)}
\int_{r\vert\bar{\sigma}_n^{(r)}\vert}^{n\vert\bar{\sigma}_n^{(r)}\vert}
t^lL\Big(\frac{t}{\vert\bar{\sigma}_n^{(r)}\vert}\Big)
\exp\Big(-tsgn(\bar{\sigma}_n^{(r)}\Big)dt,\quad
  l>0, \quad n\to \infty \la{211}\end{equation}

\nin correspondingly. Next,  in (\ref{kar}) we set
$z=e^{-\vert\sigma\vert },$ so that $1-z\sim \vert\sigma\vert,$ as $\sigma\to
0,$ and apply the  integral test to the sum $S(z)$ in the LHS, to
obtain \be \int_{ \vert\sigma\vert}^\infty
t^{l}L\Big(\frac{t}{\vert\sigma\vert}\Big)e^{-t}dt\sim \Gamma(l+1)
L\Big(\frac{1}{\vert\sigma\vert}\Big), \quad l>0, \quad \sigma\to
0. \la{h0}\end{equation}

\nin Now we are in a position to establish  asymptotic formulae for the three key parameters $\sigma, B^2 $ and
$\rho $ of the problem considered. To facilitate the understanding of the forthcoming asymptotic formulae  we
 make the following

\nin \begin{rem}\label{09} Combining  (\ref{21}) and (\ref{211})
 with (\ref{h0}), it
is easy to see that, for all $\beta,\quad 0<\beta<1$,
  both $\underline{\sigma}_n^{(r)},\bar{\sigma}_n^{(r)}$
are \be\le \tilde{ L }(n)n^{-\frac{1}{l+1}}, \quad l>0, \la{1000}\end{equation}

\nin as $n\to \infty,$ where $\tilde{ L }$ is a s.v. function induced by the given s.v.
 function $L.$ We will show in due course that $\frac{1}{l+1}$ is
 a threshold value  in the context of the problem
 considered.
\end{rem}
\vskip .5cm

 \nin It is plain that our objective requires the
derivation of asymptotic formulae for the integrals in the RHS's
of (\ref{21}), (\ref{211}).  The fact that $\sigma$ depends  on
$n$ does not allow the straightforward application of \refm[kar].
 To achieve the above goal
we
 make  use of the following fundamental fact in
the theory of s.v. functions.

\begin{prop}\label{03} (\cite{bing},Theorem 1.5.2, p. 22.) \nin For
any $b>0$ and any s.v. $L,$ the convergence

\be \Phi(x, \lambda):= \frac {(x\lambda)^{\delta} L(\lambda x)}{x^{\delta} L(x)}- \lambda^{\delta}\to 0, \quad
{as} \quad x\to \infty \la{1001}
\end{equation}

\nin is uniform in $\lambda \in [b, \infty), $ if $\delta<0,$  and is uniform in $\lambda \in (0, b],$  if
$\delta>0$ and if the function $x^{-\delta}L(x)$ is   locally bounded on $[0,\infty).$
\end{prop}
\nin Based on this result we prove now the following Abelian summation theorem which is a  version of Proposition
4.1.2, p. 199 in \cite{bing}.

\begin{prop}\label{11}

\nin Let $0<b\le \infty$ and let $b_n\to b,\ b_n
  z_n\to \infty, \  n\to\infty.$

\nin Then

\be \int_{b_n}^{\infty} e^{-t}t^lL(tz_n)dt \sim L(b_nz_n)\int_{b_n}^{\infty} e^{-t}t^l dt, \quad l>0, \quad n\to
\infty \la{802}
\end{equation}

\nin and, assuming the function $x^{-\delta}L(x)$ is locally bounded on $[0,\infty)$ for some $\delta>0,$

\be \int_{0}^{b_n} e^{t}t^lL(tz_n)dt\sim L(b_nz_n)\int_{0}^{b_n} e^{t}t^ldt, \quad l>0, \quad n\to \infty.
\la{8020}
\end{equation}
\end{prop}
\nin {\bf Proof.} In Proposition 1 we set $x=z_nb_n,$ $\lambda=
t(b_n)^{-1}$ and write the identity \be L(\lambda
x)=\Phi(x,\lambda)\lambda^{-\delta} L(x)+L(x), \la{803}
\end{equation}
\nin where $\Phi(x,\lambda)$ is as defined in (\ref{1001}). \nin Since $\lambda\ge 1$ for all $t\ge b_n,$
Proposition 1 gives \be
 \vert\int_{b_n}^{\infty} e^{-t}t^l \Phi(z_nb_n,\frac{t}{b_n})
\Big(\frac{t}{b_n}\Big)^{-\delta} dt\vert \le \epsilon b_n^{\delta}\int_{b_n}^{\infty} e^{-t}t^{l-\delta}dt,
 \quad l>0,
\la{804}
\end{equation}
\nin for all $\epsilon>0, \delta <0$ and all sufficiently large
$n.$ In (\ref{804}) we have \be
b_n^{\delta}\int_{b_n}^{\infty}e^{-t}t^{l-\delta}dt\sim \left
\{\begin{array}{ll} h, & {\rm if~} b<\infty\cr
e^{-b_n}b_n^{l+\delta} , & {\rm if~}  b=\infty,
\end{array}
\right. \la{805}
\end{equation}
\nin where $h=b^{\delta}\int_{b}^{\infty}e^{-t}t^{l-\delta}dt<\infty.$ Hence, the RHS of \refm[804] tends to $0$,
as $n\to \infty.$ Now we substitute (\ref{803}) into the LHS of (\ref{802}) to get the first assertion. The
assertion (\ref{8020}) is proved in the same manner, by applying Proposition 1 in the case $\delta>0.$ \qed

\nin  From now on, we set $\underline{r}=n^{\underline\beta},\
\bar{r}=n^{\bar{\beta}}$  and assume that the limit $d: =\lim_{n\to \infty} L(n),\ 0\le d \le \infty$ exists.
The forthcoming  assertions tell us that the latter assumption  plays a role only for the description of
the behavior of the model at the critical point.
In the case when the limit does not exist, the above description
 can be obtained  in terms of partial limits of $L(n)$, as $n\to
 \infty$.

Proposition 2 will be repeatedly used for derivation of asymptotic
formulae for the key parameters.

\begin{lemma}\label{ty}

\nin (a) Let $(l+1)^{-1}<\underline{ \beta}\le 1$ and $0\le
\bar{\beta}<(l+1)^{-1}. $ Then \be
\underline{\sigma}^{(\underline{r})}_n\sim\bar{\sigma}^{(\bar{r})}_n)\sim
\Big(\Gamma(l+1)\Big)^{\frac{1}{l+1}}n^{-\frac{1}{l+1}}L_1(n),\quad
l>0, \quad n\to \infty, \la{rt} \end{equation} \nin where
 $L_1$ is a s.v. function determined  by the s.v. function $L$ via the
relationship \be \frac{1}{L_1(n^{l+1})}\sim \Big(L^{\frac{1}{l+1}}(n)\Big)^*,\quad \quad n\to \infty. \la{*}
\end{equation}
\nin (b) Let the function $x^{-\delta}L(x), \ \delta>0$ be locally
bounded on $[0,\infty)$ and let $ 0<\underline{\beta}<(l+1)^{-1}$
and
 $(l+1)^{-1}<\bar{\beta}<1.$ Then
\be \underline{\sigma}^{(\underline{r})}_n\sim -\frac{\underline\gamma \log n}
{n^{\underline\beta}}(1+\underline{\delta}_n), \quad l>0, \quad n\to \infty, \la{92}
\end{equation}

\nin where \be \underline{\gamma}=\underline{\gamma}_n=1-(l+1)\underline{\beta}-\frac{\log L(\underline{r})}{\log
n},\quad\quad \underline{\delta}_n=\frac{\log(\underline{\gamma}\log n)}{\underline{\gamma}\log n}, \la{93}
\end{equation}

\nin while


\be \bar{\sigma}^{(\bar{r})}_n\sim \frac{\bar{\gamma}\log n}{n^{\bar{\beta}}}(1-\bar{\delta_n}), \quad l>0, \quad
n\to \infty, \la{921}
\end{equation}
\nin where \be \bar{\gamma}=\bar{\gamma}_n= (l+1)\bar{\beta}-1+\frac{\log L(\bar{r})}{\log n},\quad \quad
\bar{\delta}_n=\frac{\log(\bar{\gamma}\log n)}{\bar{\gamma}\log n}. \la{931}
\end{equation}

\nin (c) Let $ \un{\beta}=\bar{\beta}= (l+1)^{-1}.$ Then the following three cases should be distinguished:

\nin (i) If  $ 0< d<\infty,$ then \be
\underline\sigma^{(\un{r})}_n\sim
\underline{A}n^{-\frac{1}{l+1}}L_1(n), \quad l>0, \quad n\to
\infty \la{rt1}
\end{equation}
\nin and \be \bar{\sigma}^{(\bar{r})}_n) \sim
\bar{A}n^{-\frac{1}{l+1}}L_1(n),
 \quad l>0, \quad n\to \infty,
 \la{rt11} \end{equation}

 \nin
  where  $L_1$ is a
s.v. function given by (\ref{*}), while $\un{A}, \bar{A}>0$ are
the unique solutions of the equations \be
 \underline{A}^{l+1}=
\int_0^{\underline{A}d^{\frac{1}{l+1}}} t^le^{-t}dt \la{910}\end{equation}

\nin and  \be \bar{A}^{l+1}=\int_{\bar{A}d^{\frac{1}{l+1}}}^{\infty} t^le^{-t}dt \la{911}
\end{equation}

\nin correspondingly.

\nin (ii) If $d= 0,$ then $\bar{\sigma}_n^{(\bar{r})}$ is given by
(\ref{rt}), while $\underline{\sigma}_n^{(\underline{r})}$ is
given by (\ref{92}), (\ref{93}).

\nin (iii) If $d= \infty, $
 then
$\bar{\sigma}_n^{(\bar{r})}$ is given  by (\ref{921}), (\ref{931}), while $\underline{\sigma}_n^{(\underline{r})}$
is given by (\ref{rt}),(\ref{*}).
\end{lemma}

 \nin {\bf Proof.}  Since the equations
 (\ref{sig}), (\ref{sig1}) have unique solutions (in $\sigma$),
 it  suffices to check that the stated asymptotic formulae
 satisfy (\ref{21}), (\ref{211}).\qed

\begin{cor}\label{yu}
 Let $\un{\sigma}_n^{(\un{r})}, \bar{\sigma}_n^{(\bar{r})}$
be given as in Lemma 2. Then, as $n\to \infty$,

\be \un{B}^2\sim h
  \begin{cases}
  n\left(\un{\sigma}_n^{(\un{r})}\right)^{-1}, & \text{if}\ \ (l+1)^{-1}<\un\beta\le 1
\ \ \  \text{or}\ \ \un{\beta}=(l+1)^{-1}, \ d\neq 0 \\
  n\un{r}, & \text{if } \ 0<\un{\beta}<(l+1)^{-1}\ \ \text{or}\ \ \un{ \beta}=(l+1)^{-1}, \ d=
  0,
  \end{cases}
   \la{811}
\end{equation}

\be \bar{B}^2\sim h
  \begin{cases}
  n\left(\bar{\sigma}_n^{(\bar{r})}\right)^{-1},   & \text{if}\
 \  0\le\bar{\beta}<(l+1)^{-1}
\  \ \ \text{or}\   \ \bar {\beta}=(l+1)^{-1}, \ d\neq \infty\\
   n\bar{r}, & \text{if } \ \ (l+1)^{-1}<\bar{\beta}<1\ \ \
  \text{or} \ \ \bar{\beta}=(l+1)^{-1}, \ d= \infty,
   \end{cases}
   \la{812}
\end{equation}

\be \un{\rho}\sim h\frac{(\un{B}^2)^2}{n}, \quad \bar{\rho}\sim
h\frac{(\bar{B}^2)^2}{n}, \quad \quad n\to \infty.
\la{813}\end{equation}
\end{cor}

\begin{rem}\label{qs} Lemma 2 and Corollary 1 show that
$\beta=\frac{1}{l+1}$ is the critical value for the three key
parameters  $\sigma, B^2$ and $\rho$. We will see later on that
this fact has a crucial  influence on  the asymptotic behavior of
$\un{c}_n^{(r)}$ and $\bar{c}_n^{(r)}$.
\end{rem}
\vskip.5cm
 \nin Corollary 1 implies that the following
 weaker  (the third moment
 $\rho=\sum_{k=1}^n (X_k-EX_k)^3<\sum_{k=1}^n \vert X_k-EX_k\vert^3,$ see \cite{frgr},
 p.278),
 form of Lyapunov's condition (see \cite{fel}, p.278)
  holds for the sums
$\un{Y}_n^{(\un{r})}$ and $\bar{Y}_n^{(\bar{r})}$ of random variables defined by (\ref{pr}):

\be \frac{\un{\rho}}{\un{B}^3}\to 0, \quad \frac{\bar{\rho}}{\bar{B}^3}\to 0, \quad n\to \infty. \la{814*}
\end{equation}

\nin Recall that Lyapunov's condition  is sufficient for the convergence to the normal law in the central limit
theorem for independent random variables. Our next result shows that, for the triangular array considered,
 even a weaker form \refm[814*] of this condition is
sufficient for the same convergence in the local limit theorem. \vskip .5cm

\begin{thm}\label{3r} {\bf : Local limit theorem.}

\nin Let $a\in{\cal F}(l),  l> 0,$ and let $\un{\sigma}_n^{(\un{r})},\ \bar{\sigma}_n^{(\bar r)}$ be as in Lemma 2.
Then \be Pr(\un{Y}_{n}^{(\un{r})}=n)\sim (2\pi \un{B}^2)^{-\frac{1}{2}}, \quad n\to \infty,
\end{equation}
\nin

\be Pr(\bar{Y}_{n}^{(\bar{r})}=n)\sim (2\pi
\bar{B}^2)^{-\frac{1}{2}}, \quad n\to \infty.
\end{equation}
\end{thm}
\nin {\bf Proof.} Our objective will be to derive the asymptotic behavior of the integrals $\un{T}$ and $\bar{T}$
given by (\ref{333}) and (\ref{3333}) respectively. The integrands in (\ref{333}) and (\ref{3333}) are periodic
with  period $1.$ So, for any $\alpha_0$, $0<\alpha_0\le 1/2,$ the integrals can be written as

\be \un{T}=\un{T}_1 +\un{ T}_2, \quad \bar{T}=\bar{T}_1 + \bar{T}_2, \label{8111}
\end{equation}

\nin where $\un{T}_1=\un{T}_1(\alpha_0)$,
$\bar{T}_1=\bar{T}_1(\alpha_0)$ and $ T_2=T_2(\alpha_0),$
$\bar{T}_2=\bar{T}_2(\alpha_0)$ are integrals of the integrands in
(\ref{333}), (\ref{3333}) over the sets $[-\alpha_0,\alpha_0]$
and $[-1/2,-\alpha_0]\cup [\alpha_0,1/2]$ respectively. Following
the approach of \cite{frp}, \cite{frgr}, we will first show that
for an appropriate choice of $\alpha_0 =\alpha_0(n),$ the main
contributions, as $n\to \infty$, to $\un{T}$ and  $\bar{T}$ come
from $\un{T}_1,$ and $\bar{T}_1 $ respectively. From (\ref{p})-
(\ref{v2}) we have for $\alpha\in R,$

\be \varphi_j(\alpha)= \sum_{k=0}^{\infty} \frac{\left(a_j
e^{-j\sigma}e^{2\pi i\alpha j}\right )^k}
{k!\exp\left(a_je^{-j\sigma}\right)}=
\exp\left(a_je^{-j\sigma}\left(e^{2\pi i\alpha j}-1\right
)\right)\label{812*} \end{equation}

\nin  and \be \un{\varphi}^{(\un{r})} (\alpha)=
\exp\left(\sum_{j=1}^{\un{r}}
a_je^{-j\un{\sigma}_n^{(\un{r})}}\left(e^{2\pi i\alpha j} -1\right
)\right),
 \label{8121} \end{equation}

\be \bar{\varphi}^{(\bar{r})} (\alpha)=
\exp\left(\sum_{j=\bar{r}}^n
a_je^{-j\bar{\sigma}_n^{(\bar{r})}}\left(e^{2\pi i\alpha j}
-1\right )\right). \label{813*} \end{equation}

\nin Substituting    the Taylor expansion (in $\alpha$)

\be e^{2\pi i\alpha j}-1= 2\pi i\alpha j-2\pi^2\alpha^2j^2 +O(\alpha^3j^3), \quad \mbox{as $\alpha\to 0$},
\label{816}
\end{equation}

\nin in (\ref{8121}), (\ref{813*}) and employing  (\ref{sig}),
(\ref{sig1}),
 gives \be \un{\varphi}^{(\un{r})}(\alpha)e^{-2\pi i\alpha n} =
\exp\left(-2\pi^2\alpha^2\un{B}^2 +
 O(\alpha^3\un{\rho})\right), \quad
\mbox{as \ $\alpha\to 0$}, \label{814} \end{equation}

\be \bar{\varphi}^{(\bar{r})}(\alpha)e^{-2\pi i\alpha n} = \exp\left(-2\pi^2\alpha^2\bar{B}^2 +
 O(\alpha^3\bar{\rho})\right), \quad
\mbox{as \  $\alpha\to 0$}. \label{8141} \end{equation}
\nin We write now \be \un\alpha_0^3
\un{\rho}=(\un\alpha_0\un{B})^3\frac{\un{\rho}}{\un{B}^3}, \quad \bar{\alpha}_0^3
\bar{\rho}=(\bar{\alpha}_0\bar{B})^3 \frac{\bar{\rho}}{\bar{B}^3}
\end{equation}
\nin to conclude that, by (\ref{814*}), there exist
$\un{\alpha}_0=\un{\alpha}_0(n), \
\bar{\alpha}_0=\bar{\alpha}_0(n)$ such that \be \lim _{n\to\infty}
\un{\alpha}_0 \un B=\lim_{n\to \infty} \bar{\alpha}_0\bar{B}
=+\infty \label{815} \end{equation} \nin and \be \lim_{n\to\infty}
\un{\alpha}_0^3 \un{\rho}
=\lim_{n\to\infty}\bar{\alpha}_0^3\bar{\rho} =0 .\label{816*}
\end{equation} \nin We see from (\ref{816*}) that $ \un{\alpha}_0,\
\bar{\alpha}_0\to 0, \ n\to \infty,$ because $\un{\rho},\
\bar{\rho}\to \infty, \ n\to \infty,$ by \refm[813]. Also note the
fact that (\ref{816}) holds for all $\alpha\in[-\un{\alpha}_0,\un
{\alpha}_0] \bigcup[-\bar{\alpha}_0,\bar{\alpha}_0].$ \nin As a
result, we arrive at  the  asymptotic formulae
 for the integrals
$\un{T}_1, \bar{T}_1:$

 \ber \un{T}_1\sim
\int_{-\un{\alpha}_0}^{\un{\alpha}_0}
\exp\left(-2\pi^2\alpha^2\un{B}^2\right)d\alpha=\nonumber\\
\frac{1}{2\pi\un{ B}} \int_{-2\pi\alpha_0\un{B}}^{2\pi
\alpha_0\un{B}} \exp(-\frac{z^2}{2}) dz \sim \frac{1}{\sqrt{2\pi
\un{B}^2}},  \quad n\to \infty ,\label{817}
\end{eqnarray}

\ber \bar{T}_1\sim \int_{-\bar{\alpha}_0}^{\bar{\alpha}_0}
\exp\left(-2\pi^2\alpha^2\bar{B}^2\right)d\alpha=\nonumber\\
\frac{1}{2\pi \bar{B}} \int_{-2\pi\bar{\alpha}_0{\bar{B}}}^{2\pi
\bar{\alpha}_0\bar{B}} \exp(-\frac{z^2}{2}) dz \sim
\frac{1}{\sqrt{2\pi \bar{B}^2}},  \quad n\to \infty .\label{818}
\end{eqnarray} \nin Now  we turn to the estimation, as $ n\to
\infty,$ of the integrals $\un{T}_2$, $\bar{T}_2.$

\nin We have

\be \vert \un{T}_2\vert =2\vert \int_{\un{ \alpha}_0}^{
1/2}\un{\varphi}^{(\un r)}(\alpha) e^{-2\pi i\alpha
n}d\alpha\vert, \quad \vert\bar{T}_2\vert =2\vert \int_{\bar{
\alpha}_0}^{ 1/2}\bar{\varphi}^{(\bar r)}(\alpha) e^{-2\pi i\alpha
n}d\alpha\vert. \la{820}\end{equation}

\nin It follows from (\ref{8121}), (\ref{813*}) that

\be \vert \un{\varphi}^{(\un{r})}
(\alpha)\vert=\exp\left(-2\sum_{j=1}^{\un{r}}
a_je^{-j\un{\sigma}_n^{(\un{r})}}\sin^2\pi \alpha j\right), \quad
\alpha\in R, \label{821} \end{equation} \be
\vert\bar{\varphi}^{(\bar{r})} (\alpha)\vert=
\exp\left(-2\sum_{j=\bar{r}}^n
a_je^{-j\bar{\sigma}_n^{(\bar{r})}}\sin^2\pi \alpha j\right),
\quad \alpha\in R. \label{822} \end{equation}
\nin We denote \be \un{V}_n^{(\un{r})}(\alpha)=2\sum_{j=1}^{\un{r}} a_je^{-j\sigma_n^{(\un{r})}}\sin^2\pi\alpha j,
\quad \un\alpha_0\le \alpha\le 1/2, \la{823}
\end{equation}

\be \bar{V}_n^{(\bar{r})}(\alpha)=2\sum_{j=\bar{r}}^n a_je^{-j\bar{\sigma}_n^{(\bar r)}}\sin^2\pi\alpha j, \quad
\bar{\alpha}_0\le \alpha\le 1/2. \la{823*}
\end{equation}

\nin For the sake of estimating the sums $\un{V}_n^{(\un{r} )}, \
\bar{V}_n^{(\bar r)}$ we make use of the following inequality
proven in \cite{frp}: \be 2\sum_{j=p}^{p+k-1} sin^2 \pi \alpha
j\ge \frac{k}{2} \min\{1,(\alpha k)^2\}, \quad \vert
\alpha\vert\le 1/2,\quad \forall k\ge 2,\ p\ge 1. \la{824}
\end{equation}

\nin
  We set
 \be \un{\alpha}_0^2=\frac{\log^4(\un{ B}^2)}{\un{B}^2}, \quad \bar{\alpha}_0^2=
 \frac{\log^4 (\bar{B}^2)}{\bar{B}^2},  \la{830}\end{equation} and
apply (\ref{824}) with \be \un{k}=l\left(\vert\un{\sigma}_n^{(\un{r})}\vert\right)^{-1}, \quad
\bar{k}=l\left(\vert\bar{\sigma}_n^{(\bar{r})}\vert\right)^{-1} \la{828}
\end{equation}
\nin and different $\un{p},\bar{p}.$ \nin (Note that under the
choice (\ref{830}) of $\un{\alpha}_0, \bar{\alpha}_0,$ the
conditions (\ref{815}), (\ref{816*}) indeed hold.) Treating
separately the cases (a), (b) and (c) in Lemma 2, we are able to
show that

 $$ e^{-\un{V}_n^{(\un{r})}(\alpha)}=o(\un{B}^{-1}), \quad
\un{\alpha}_0\le \alpha\le 1/2, \quad n\to \infty \quad and$$ \be
e^{-\bar{V}_n^{(\bar{r})}(\alpha)}=o(\bar{B}^{-1}), \quad
\bar{\alpha_0}\le \alpha\le 1/2, \quad n\to \infty.
\la{8311}\end{equation}
\qed

\nin \begin{cor}\label{u8}{\bf :Asymptotic formulae for
$\un{c}_n^{(\un{r})},
 \bar{c}_n^{(\bar{r})}.$}

\nin Let $\un{r}=n^{\un{\beta}}, \ 0<\un{\beta}\le 1,$ and $\bar{r}=n^{\bar{\beta}}, \ 0\le \bar{\beta}< 1.$ Then

\be \un{c}_n^{(\un{r})}\sim (2\pi \un{B}^2
)^{-\frac{1}{2}}\exp\Big(\un{S}_n^{(\un{r})}(e^{-\un{\sigma}_n^{(\un{r})}}) +n\un{\sigma}_n^{(\un{r})}\Big) ,
\quad n\to \infty, \la{836}
\end{equation}

\be \bar{c}_n^{(\bar{r})} \sim (2\pi \bar{B}^2)^{-\frac{1}{2}}
\exp\Big(\bar{S}_n^{(\bar{r})}(e^{-\bar{\sigma}_n^{(\bar{r})}})+ n\bar{\sigma}_n^{(\bar{r})}\Big), \quad n\to
\infty, \la{837}
\end{equation}
\nin where $$ \un{S}_n^{(\un{r})}(e^{-\un{\sigma}_n^{(\un{r})}})\sim h \frac{n^2}{\un{B}^2}, $$ \be
\bar{S}_n^{(\bar{r})}(e^{-\bar{\sigma}_n^{(\bar{r})}})
\sim h\frac{n^2}{\bar{B}^2}, \quad n\to \infty. \la{838}
\end{equation}
\end{cor}

\nin {\bf Proof.} By Theorem 1 and  (\ref{3y}), (\ref{3y1}) we get
 the asymptotic expressions (\ref{836}), (\ref{837}), while
(\ref{838}) is obtained with the help of the integral test,
Proposition 2, Lemma 2 and Corollary 1. \qed

\vskip .5cm


\begin{thm}\label{4r}{\bf :The limiting behavior of
$\un{d}_n^{(\un{r})}, \ \bar{d}_n^{(\bar r)}.$ }

\nin Denote $\un{c}_n^{(n)}=c_n.$

\nin (i)
 Let $ \un{r}=n^{\un{\beta}},\quad 0\le \un{\beta}\le 1.$ Then

\nin \be \lim_{n\to \infty}\un{d}_n^{(\un{r})} =
  \begin{cases}
    0, & \text{if}\quad 0\le \un{\beta}<(l+1)^{-1} \
    \ \text{or}\ \ \un{\beta}=(l+1)^{-1}, \ d<\infty\\
    1, & \text{if} \quad (l+1)^{-1} < \un{\beta}\le 1 \ \ \text{or}\
    \un{\beta}=(l+1)^{-1}, \ d=\infty.
  \end{cases}
\la{859}
\end{equation}

\nin (ii) Let  $\bar{r}\ge 2.$ Then \be \lim_{n\to
\infty}\bar{d}_n^{(\bar r)} =
  \begin{cases}
   0, & \text{if}\quad \bar{r}=n^{\bar{\beta}}, \quad 0<\bar{ \beta}\le 1\\
    \exp\Big(-\sum_{j=1}^{\bar{r}-1} a_j\Big),  & \text{if}
    \quad \bar{r}\ge 2 \quad \text{is a fixed
    number}.
  \end{cases}
  \la{8591}
\end{equation}
\end{thm}

\nin {\bf Proof:}
 (i)   Denote
\be \Delta_n^{(\un{r})}= \un{S}_n^{(\un{r})}\left(e^{-\un{\sigma}_n^{(\un{r})}}\right)-
\un{S}_n^{(n)}\left(e^{-\un{\sigma}_n^{(n)}}\right) +
 n\left(\un{\sigma}_n^{(\un{r})}-\un{\sigma}_n^{(n)}\right).
\end{equation}
\nin Our objective will be to demonstrate that

\be \lim_{n\to \infty} \Delta_n^{(\un{r})}=
  \begin{cases}
    -\infty, & \text{if}\quad   0<\un{ \beta}<(l+1)^{-1} \quad \text{or}\
    \ \un{\beta}=(l+1)^{-1}, \ d<\infty\\
    0, & \text{if}\quad (l+1)^{-1}<\un{\beta}\le 1\quad \text{or}\
   \ \un{\beta}=(l+1)^{-1}, \ d=\infty.
  \end{cases}
  \la{8391}
\end{equation}

\nin The first part of (\ref{8391}) follows in a straightforward
way  from the preceding asymptotic analysis. In the case
corresponding to the second part, \be \un{\sigma}_n^{(n)}\sim
\un{\sigma}_n^{(\un{r})}, \quad n\to \infty, \la{8411}
\end{equation}
and, consequently,

  \be
\un{S}_n^{(\un{r})}\left(e^{-\un{\sigma}_n^{(\un{r})}}\right)\sim
\un{S}_n^{(n)}\left(e^{-\un{\sigma}_n^{(n)}}\right), \quad n\to
\infty. \la{850}
\end{equation}
\nin Therefore, here a more subtle analysis is required. From the
identity \be\un{ M}_n^{(\un{r})}(\un{\sigma}_n^{(\un{r})})-\un{
M}_n^{(n)}(\un{\sigma}_n^{(n)})=0, \ \ r=n^\beta,\
1\ge\beta\ge\frac{1}{l+1},\ \ n=2,3,\ldots \end{equation}

\nin and the fact that
$\un{\sigma}_n^{(n)}>\un{\sigma}_n^{(\un{r})}\ge 0,$ we derive
that $n(\un{\sigma}_n^{(n)}-\un{\sigma}_n^{(\un{r})})\to 0, \ n\to
\infty. $ Similarly, it can be proven that in the case considered
\be \un{S}_n^{(\un{r})}\left(e^{-\un{\sigma}_n^{(\un{r})}}\right)-
\un{S}_n^{(n)}\left(e^{-\un{\sigma}_n^{(n)}}\right)\to 0,\  n\to
\infty.
\end{equation}

\nin Combining this with the asymptotic formulae in Corollary 1 and Corollary 2, proves the
second part of  (\ref{8391}).

\nin (ii) We outline only  the proof of the second part of
\refm[8591].
\nin  We have  for a fixed $\bar r,$
\be
\bar{\sigma}_n^{(\bar{r})}\sim \bar{\sigma}_n^{(1)}, \quad n\to
\infty,\quad and \quad \lim_{n\to \infty} n\Big(
\bar{\sigma}_n^{(\bar r)}-\bar{\sigma}_n^{(1)}\Big)= 0.\la{nm}
\end{equation}
\nin Next we write \be \bar{S}_n^{(\bar
r)}\left(e^{-\bar{\sigma}_n^{(\bar r)}}\right)-
\un{S}_n^{(n)}\left(e^{-\bar{\sigma}_n^{(1)}}\right)= \sum_{j=1}^n
a_j e^{-\bar{\sigma}_n^{(1)} j}\left( e^{-(\bar{\sigma}_n^{(\bar r
)}-\bar{\sigma}_n^{(1)}) j}-1\right) -\sum_{j=1}^{\bar{r}-1} a_j
e^{-\bar{\sigma}_n^{(\bar r)}j}. \la{wc}
\end{equation}
\nin Since $\bar{\sigma}_n^{(1)}>\bar{\sigma}_n^{(\bar{r})}\ge 0$
and
\be e^{-(\bar{\sigma}_n^{(1)}-\bar{\sigma}_n^{(\bar r)})j} -1= \Big(\bar{\sigma}_n^{(1)}-\bar{\sigma}_n^{(\bar
r)}\Big)j(1-\delta_n), \quad 1\le j\le n, \la{we}
\end{equation}

\nin where $\delta_n=\delta_n(j)\to 0, \quad n\to \infty,$
uniformly in $1\le j\le n,$ we get the desired claim. \qed

\nin \section{Application 1: Reversible coagulation-fragmentation processes.}

\setcounter{equation}{0}

\nin We follow the formulation of the model given in  \cite{dgg}. A population of $n$ particles is partitioned
into groups of various sizes that undergo stochastic  evolutions (in time) of coagulation  and
fragmentation. There are only two  possible interactions:
  coagulation of two groups into
one, and fragmentation of one group into two groups. The  process of
coagulation-fragmentation(CFP) is a time-homogeneous interacting
particle system  $\varphi_t, \ t\ge 0,$  defined as follows.
 For a given $n,$ denote by $\eta=( k_1, \ldots, k_n)$  a
partition of the whole population $n$ into $k_i$ groups of size
$i, \ i=1, 2, \ldots,n,$ where the numbers of groups $k_i \ge 0$
are subject to the condition:

\be \sum_{i=1}^{n} i k_i = n, \la{c00}
\end{equation}

\nin called the total mass conservation law. The finite set $\Omega_n = \{\eta \}$ of all partitions of $n$ is the
state space of the process $\varphi_t, t\ge 0$. The rates of the infinitesimal (in time)
 transitions (= flips)
 are assumed to depend
only on the sizes of the interacting groups, and are given by two functions $\psi$ and $\phi:$

\begin{enumerate}
\item
For $i$ and $j$ such that $i+j \le n,$ the rate of coagulation, $(i,
j) \to (i+j),$ of two groups of sizes $i$ and $j$ into one group
of size $i+j, $   equals  $\psi(i, j)$.

\item
The rate of fragmentation, $(i+j) \to (i, j)$, of a group of size
$i+j$ into two groups of sizes $i$ and $j,$ equals  $\phi(i,
j)$.
\end{enumerate}

\nin Hereafter, we refer to the coagulation and fragmentation
rates $\psi$ and $\phi$ as intensities. The intensities are
required to satisfy $\psi(i, j) = \psi(j, i) \ge 0$ and $\phi(i,
j) = \phi(j, i) \ge 0$. We also make the natural assumption
 that the total intensities of merging $\Psi(i,j; \eta)$ and splitting
$\Phi(i,j;\eta)$ at a configuration $\eta\in \Omega_n$ are given
by:

\ber \Psi(i, j;\eta)= \Psi(i, j; k_i, k_j) = \psi(i, j) ~ (k_i
k_j)^\gamma, \quad i \ne j,  \quad 2\le i+j \le n, \nonumber \\
\Psi(i, i;\eta)= \Psi(i, i; k_i, k_i) =\psi(i,i)
(k_i(k_i-1))^\gamma, \quad 2\le 2i\le n, \nonumber \\ \Phi(i,
j;\eta)= \Phi(i, j; k_i, k_j) = \phi(i, j) ~ (k_{i+j})^\gamma,
\quad 2\le i+j \le n, \nonumber \\ \label {formula_fl}
\end{eqnarray}
\nin where  $\gamma >0.$  Note that the  case $\gamma = 1$ corresponds to the mass action kinetics.

\nin In this paper, we study only reversible CFP's with nonzero
intensities. It is  known ( \cite{dgg}, \cite{kel}) that the class
of such processes is characterized by the following property of
the ratio of their intensities:

\be \frac{\psi(i,j)}{\phi(i,j)}= \frac{a_{i+j}}{a_ia_j}, \quad 2\le i+j \le n, \la{c1}
\end{equation}

\nin where $a=\{a_i=a(i)\}, \ i=1,2,\ldots$ is a positive function. It is also known

\nin (\cite{dgg}, \cite{kel}) that, under (\ref{c1}), the unique
invariant measure $\mu_n$ on $\Omega_n$ is given by

\ber \mu_n (\eta) = C_n \frac{a_{1}^{k_{1}} a_{2}^{k_{2}} \ldots a_{n}^{k_n}} {(k_{1}!)^\gamma (k_{2}!)^\gamma ,
\ldots (k_{n}!)^\gamma}, \quad \eta=(k_1,\ldots,k_n)\in \Omega_n.
\non \\
\label{c2}
\end{eqnarray}
 \nin
Here $C_n^{-1}(a)=C_n^{-1}:=c_n$ is the partition function
  for the probability measure $\mu_n, \ n\ge 1$:

\ber c_n = \sum_{\eta\in \Omega_n} \frac{a_1^{k_1}
a_2^{k_2} \ldots a_n^{k_n}} {(k_1!k_2!  \ldots k_n!)^\gamma},
\quad \eta=(k_1,\ldots,k_n)\in \Omega_n. \label{c3} \end{eqnarray}
\nin The measure $\mu_n$ is the steady state of the reversible CFP
considered. So, \refm[c2] tells us that for a fixed $n,$ the
steady state is determined  by $ n$ values of the function $a.$ In
view of this, it is natural to call $a$ the parameter function of
the process. Note that in contrast to the above, the transient
behavior of the CFP's considered depends on the intensities
$\psi$ and $\phi$, rather than on their ratios.

\begin{rem}\label{1d} The measure $\mu_n$ is invariant under the
following transformation of the parameter function $a$.  Define
the family of operators $H_h, \ h>0$ on a set of parameter
functions $a:$ $$ (H_h a)(j)=h^j a_j,\quad  j=1,2,\ldots, \quad
h>0.$$ It follows from \refm[c00] and \refm[c2] that (with the
obvious abuse of notation) \be H_h\mu_n=\mu_n, \quad
h>0.\la{c900}\end{equation} This says that  all results of the
present paper are extended to  the
 class of parameter functions
 $\{H_h a: h>0,\ a\in{\cal F}(l), \ l> 0\}. $
\refm[c900] also explains the possibility of introducing a free parameter for the treatment of problems related to
$\mu_n .$
\end{rem}
\vskip .5cm \nin Our study is devoted exclusively to the steady state of the above CFP's, in the case when in
\refm[c2], $\gamma=1$ and $n\to \infty.$
 Treating $ S$ given by \refm[gen] as a generating function of the positive sequence $\{a_n\}_1^\infty,$
 such that
the radius of convergence of the series \refm[gen] equals $1$, it is known (see e.g. \cite{dgg}) that
 the sequence $\{c_n\}_0^\infty$ in \refm[c3] is generated by the function
 $g$ defined by
\be g(z) = e^{S(z)} = \sum_{n=0}^{\infty} c_n z^n,\quad \vert z\vert< 1.\la{c4}
\end{equation}
\nin To formulate the problem of clustering in the  setting of CFP, we
 define on  the  probability space $(\Omega_n,\mu_n)$
the random variables  $K_i=K_i^{(n)}(\eta) =$
 the number of groups of size $i, \ i=1,\ldots, n,$ in a random
partition $\eta\in\Omega_n,$ and let $\bar{q}_n=\bar{q}_n(\eta),$ (resp. $\un q_n= \underline{q}_n(\eta)$) be
 the size of the largest (resp. smallest) group.
  We will be interested in the probabilities
$Pr(\bar{q}_n\le r)$ and $Pr(\underline{q}_n\ge  r).$ Making use
of the notation in \refm[gen2], we have $$ \un{c}^{(
r)}_n=\sum_{\eta\in \Omega_n:\bar{q}_n(\eta)\le   r}
\frac{a_1^{k_1} a_2^{k_2} \ldots a_n^{k_n}} {k_1!k_2!  \ldots
k_n!},$$ \be \bar{c}^{(r)}_n=\sum_{\eta\in
\Omega_n:\underline{q}_n(\eta)\ge r} \frac{a_1^{k_1} a_2^{k_2}
\ldots a_n^{k_n}} {k_1!k_2!  \ldots k_n!}. \la{c5}
\end{equation}

\nin This gives

\be Pr(\bar{q}_n\le r)=\frac{\un c^{(r)}_n}{c_n}, \quad
Pr(\underline{q}_n\ge r))= \frac{\bar{c}^{(r)}_n}{c_n}. \la{c6}
\end{equation}

\nin We  assume now that  $a\in  {\cal F}(l), \  l> 0,$ and
$r=n^\beta, \ 0\le \beta \le 1.$ Then Theorem 3 admits the
following interpretation :

\nin  \be \lim_{n\to \infty} Pr(\bar{q}_n\le n^\beta) =
  \begin{cases}
    0, & \text{if}\quad 0\le \beta<(l+1)^{-1} \quad \text{or}
    \quad \beta=(l+1)^{-1}, \ d<\infty  \\
    1, & \text{if} \quad  (l+1)^{-1}<\beta\le 1
     \quad \text{or}
    \quad \beta=(l+1)^{-1}, \ d=\infty,
  \end{cases}
\la{c7}
 \end{equation}
\nin while
  \be \lim_{n\to \infty} Pr(\underline{q}_n\ge r)=
\begin{cases}
   0, & \text{if}\quad r=n^\beta, \quad 0< \beta\le 1\\
   \exp\Big(-\sum_{j=1}^{r-1} a_j\Big) ,  & \text{if}
    \quad r\ge 2 \quad \text{is a fixed
    number}.
  \end{cases}
  \la{c71}
\end{equation}

\begin{rem}\label{1z}

\refm[c7] identifies $n^{\frac{1}{l+1}}$ as the threshold for the
limiting distribution  of the size of the largest cluster, in the
sense that  $$\frac{1}{l+1}=\inf\{ \beta: Pr(\bar{q}_n\le n^\beta
)=1\}.$$
 We
discuss the phenomenon in more details in Remark 8, in the context
of random combinatorial structures.
\end{rem}

\nin To reveal the picture of clustering  at the equilibrium of the CFP's considered, we establish one more fact.

\begin{thm}\label{3e}

\nin   Let $a\in  {\cal F}(l), \  l> 0.$ Then

\nin (i) \be \lim_{n\to \infty} Pr(n^{\frac{1}{l+1}-\epsilon} <\bar{q}_n<n^{\frac{1}{l+1}+\epsilon})=1, \quad
\forall \epsilon>0. \la{c10}
\end{equation}

\nin(ii) For all $p$ such that $n^\epsilon\le p\le n^\beta,$ with
$ \beta<\frac{1}{l+1}$  and $ \epsilon>0,$

\be \lim_{n\to \infty} Pr (K_p=0)=
  \begin{cases}
    1, & \text{if} \quad 0<l<1, \\
e^{-d},& \text{if} \quad l=1,\\
    0, & \text{if} \quad l>1.
  \end{cases}
  \la{zn}
\end{equation}

\nin (iii) For any two $s$-tuples of integers $p_1,\ldots, p_s\ge
1 $ and $k_1,\ldots, k_s\ge 0,$
 \be
\lim_{n\to \infty} Pr (K_{p_1}=k_1,\ldots,
K_{p_s}=k_s)=\prod_{j=1}^s\frac{a_{p_j}^{k_j}}{k_j!} e^{-a_{p_j}}.
\la{zv}
\end{equation}
\end{thm}
\nin {\bf Proof:} \refm[c10]   follows immediately from (\ref{c7}). Next, we have \be Pr (K_p=0)=\sum_{\eta\in
\Omega_n: k_p=0} \mu_n(\eta):=\frac{c_{n,p}}{c_n}, \quad 1\le p\le n. \la{c11}\end{equation}

\nin Denote by
 $\sigma_{n,p}, \  B_{n,p}^2,$  the key parameters
of the asymptotics of $c_{n,p}.$ Namely, $\sigma_{n,p}$ is the
unique solution in $\sigma$ of the equation \be
\un{M}^{(n)}_n(\sigma)-pa_pe^{-\sigma p}=n, \quad l>0, \la{c14}
\end{equation}
\nin and $B_{n,p}^2$ is defined correspondingly. Then in the case
$1\le p< n^\beta, \ 0< \beta<\frac{1}{l+1},$ we have
 \be
\sigma_{n,p}\sim \un{\sigma}_{n}^{(n)} \sim
\Big(\Gamma(l+1)\Big)^{\frac{1}{l+1}}n^{-\frac{1}{l+1}} L_1(n), \
\quad B_{n,p}\sim \un{B}_n^{(n)}, \quad n\to\infty, \quad
\sigma_{n,p}<\un{\sigma}_{n}^{(n)}.\la{c16}
\end{equation}
\nin  By the reasoning used for the proof of the second part of
\refm[8591], we get from \refm [c16], for $n\to \infty,$ \be 0\le
n(\sigma_{n,p}-\un{\sigma}_n^{(n)})\to 0,  \ and \
\ S_{n,p}\left(e^{-\sigma_{n,p}}\right)-
\un{S}_n^{(n)}\left(e^{-\un{\sigma}_n^{(n)}}\right)
+a_pe^{-p\sigma_{n,p}}\to 0, \ l>0,
\end{equation}
\nin which implies  \refm[zn]. Next, the relationship
 \be Pr(K_p=k_p)=\frac{a_p^{k_p}}{k_p!}Pr(K_p=0), \ 1\le p\le n,
\end{equation}
\nin implies \refm[zv] for $s=1$. For general $s$ the proof is
similar.
 \qed
\vskip .5cm
 \begin{rem} \label{2w} \refm[zv] tells us that the
random variables depicting  numbers of groups of  fixed sizes
become independent, as $n\to \infty$. This fact is in accordance
with the general principle of asymptotic independence  of
particles in mean-field models, that is commonly accepted (but not
rigorously proved) in statistical physics (see \cite{frgr} and
references therein). In the case $a\in{\cal F}(l), \ l\le 0$ the
independence principle was broadly discussed in the context of
random combinatorial structures (see  Remark 8).
\end{rem}

\nin Now we are in a position to provide  a verbal description of
the striking  feature of  clustering for large $n,$ in the case
$a\in{\cal F}(l), \ l> 0.$

 \nin
\begin{itemize}
\item With probability $1$,
there are no clusters(=groups) of sizes greater than $O(n^\beta), \ \beta>\frac{1}{l+1}$. Moreover, with
probability $1$, the size of the largest group lies  in the interval
$[n^{\frac{1}{l+1}-\epsilon},n^{\frac{1}{l+1}+\epsilon}], \ \epsilon
>0$.

\nin On the other hand:

\item  The $n$ particles are partitioned into
groups of sizes not greater than $O(n^{\frac{1}{l+1}})$ in such a way that

(i) with a positive probability there are groups of any fixed size;

(ii)  the limiting probability of having a group of a  size $p,\ p\in
[n^\epsilon,n^\beta ],$\
\ $ \epsilon>0, \ 0<\beta<\frac{1}{l+1}, $ equals
 $0,$  $1-e^{-d}$ or $1$, if \ $0<l<1,$ $l=1$ or
$l>1$ respectively.

\end{itemize}

\nin Summing up the aforementioned picture, we conclude that for
large $n$ the  distribution of clusters induced by the measure
$\mu_n$ has a threshold $n^{\frac{1}{l+1}}$.

\nin{\bf Historical remarks}

\nin It is generally accepted that the mathematical chapter of the history of CFP's  traces back to the 1917 paper
\cite{sm} by M. Smoluchowski. In this seminal work the
 mathematical theory of the process of pure coagulation of molecules
 of colloids     was proposed. A deep discussion of the physical
 context and implications of Smoluchowski model was presented in Ch.III of the
 classical work  by
 S. Chandrasekhar (1943) reprinted in \cite{wa}.
 Observe that coagulation was treated by
 Smoluchowski as a
 deterministic process.
 In the framework of this approach,  there was
  derived in \cite{sm}
   an infinite  system of
  differential equations describing the evolution in time of the
concentration of molecules of sizes $1,2,\ldots $. (Note that some authors mistakenly attribute the equations to
another paper by Smoluchowski published in 1916).

\nin Subsequently, the equations, after being generalized to allow also fragmentations of particles,
 became famous as a general
 model for processes of grouping and splitting in numerous fields

\nin  of applications. Efforts of generations of
 researchers were devoted to the intriguing mathematical
problems of existence, uniqueness and asymptotic behavior (in time and in the number of particles)  of the
solutions.

\nin It was  understood a long time ago that a stochastic context could be attributed to Smoluchowski's equations
(SE). Corresponding stochastic models were independently reintroduced, under different names, in different fields
of applications (see for details the review \cite{ald1}). The  paper \cite{ald2} attributes   to A. Marcus the
first stochastic model for pure coagulation, called the Marcus-Lushnikov process (MLP). Extensive study  of MLP was
concentrated around two subjects:
 the gelation phenomena and the relation of MLP to SE.
 (``Gelation" is the name for the
phase transition exhibited by the formation of a  giant cluster that causes the violation of the total mass
conservation law \refm[c00]).
 The main approach to these problems is based on treating the
 MLP
 as the stochastic coalescent.
 A program for investigating the relationship
between MLP (= stochastic coalescent) and SE
 was outlined by D. Aldous in \cite{ald1}.
Recent progress in this direction was made by J. Norris in \cite{nor}, who proved that under certain conditions a
sequence of stochastic coalescents converges weakly to the solution of the SE. The  theory of coalescents as a tool to
study limits of coagulation models as $n\to \infty,$ was developed by J. Pitman et. al (see e.g. \cite{ep}).
\nin Parallel to this line of research, Monte Carlo algorithms based on MLP were developed for the numerical
treatment of SE (see \cite{ ebw}) and references therein).

\nin P. Whittle \cite{wh1} proposed a reversible Markov process as
a model for Flory's theory of polymerization developed in the
1940s.
 As a result, a system of SE (in the presence
of fragmentation) was rediscovered for both deterministic and stochastic contexts (see also \cite{wh2}).
 M. Aizenmann and T. Bak \cite{a}, also motivated by
Flory's theory, proved that for the continuous (in space) version of SE with constant kernels of coagulation and
fragmentation, the free energy of the system decays exponentially as time $t \to \infty.$ This important fact
established the validity of Boltzmann's H-theorem for the time evolution of the system described by SE. Note that
a general fact of increasing
  entropy for
SE with  kernels obeying the deterministic reversibility condition was independently proven in \cite{wh1}.

 \nin The explicit
formulation of a CFP as a Markov process on the set of partitions appears in the monograph \cite{kel}, Ch.8, by
F. Kelly, which contains also the expression \refm[c2] for the equilibrium distribution of reversible CFP's in the
case $\gamma=1.$ (In \cite{kel} the model is called a clustering process.) The above formulation was reintroduced by
S. Gueron in \cite{gu} in  the context of animal grouping. As far as we know, Gueron, \cite{gu}, was the first to
notice that  SE are obtained from the Kolmogorov forward equations for the expected numbers of groups, by
neglecting correlations among the numbers  $K_p$  of groups of different sizes $p=1,2,\ldots$. R. Durrett,
B. Granovsky and S. Gueron \cite{dgg} studied the asymptotic  behavior (in n) of $EK_p$ and $cov(K_{p_1}, K_{p_2})$
at the steady state \refm[c3] with $\gamma=1$ and an arbitrary parameter function $a.$ They showed that
\be\lim_{n\to \infty} cov(K_{p_1}, K_{p_2})=0, \la{020}\end{equation} for any fixed $p_1\neq p_2,$ which agrees
(for large $n$) with the assumption of independence of group numbers of fixed sizes, in the stochastic context of
SE.  In \cite{dgg} it was also shown that for a wide class of the parameter functions $a$ and a fixed $p,$ \be
EK_p\sim k_p, \quad {as} \quad n\to \infty,\la{021}\end{equation}  where $k_p$ is the equilibrium solution of the
continuous version of SE. However, it was found that both \refm[020] and \refm[021] fail when the group size
$p=p(n)\to \infty,$ as $n\to \infty.$ The latter leads to the crucial difference in the behavior of stochastic and
deterministic
 solutions at equilibrium. It is plain that the difference between the
two models is the consequence of the mass conservation law \refm[c00] that contradicts the  assumption of
independence.

\nin In the paper \cite{frgr} by G. Freiman and B. Granovsky,
 Khintchine's probabilistic method was brought to the scenario.
With the help of this method, asymptotic formulae for the partition function  for the invariant measure \refm[c2]
were derived in the case when $a_n\sim n^{l-1}, \ l>0, \ n\to \infty.$
 In \cite{frgr} one can also find a sketch of the history of
Khintchine's method.

\nin I. Jeon \cite{j} found sufficient conditions on intensities of coagulation and fragmentation in SE under
which the  gelation phenomena occurs. Note that these conditions are not satisfied for the reversible intensities
generated by the class ${\cal F}(l), \  l> 0,$ \ of parameter functions considered in the present paper.

\nin The paper \cite{lw} by P. Laurencot and D. Wrzosek introduced a version of SE with coagulation and collisional
fragmentation. The latter means that the fragmentation ( = breakage) occurs only as a result of a collision of two
clusters.

\nin Essentially, all   stochastic and deterministic processes discussed so far  are  mean-field models, in the
sense that the rates of coagulation and fragmentation depend on the sizes of interacting groups only. J. R. Norris
\cite{nor2} formulated  a continuum version of SE in the case when the coagulation rates depend not only on the
particle masses but also on some other characteristics of the clusters (e.g., the  shape of the cluster, the types
of basic particles that form the cluster, etc).

\nin \section{Application 2: Random Combinatorial Structures(RCS).}

A combinatorial structure (CS) of a  size $n$ is defined as a union of  components (= nondecomposable elements) of
sizes $1,2,\ldots,n$,
 and by RCS we mean  the uniform
probability distribution on  the finite set of all $p_n$ CS's of size $n.$ The RCS induces  the component size
counting process ${\bf{K}}^{(n)}=(K_1^{(n)},\ldots,K_n^{(n)}),$ where $K_i=K_i^{(n)}, \ i=1,2\ldots,n$ are the
numbers of components
 (in a randomly chosen CS)
 of sizes $i=1,2,\ldots, n,$ subject to  the total mass
conservation law \refm[c00]. It was long ago understood that for  a wide class  of RCS's
 the  distribution laws $\cal{L}$ of the processes ${\bf{K}}^{(n)}$ have the
 following common feature called the conditioning relation
 (for references see the monograph \cite{abt}, Ch.2,
 by R. Arratia, A. Barbour, and
  S.Tavar\'e and \cite{arr}, by R.Arratia and S.Tavar\'e):
\be {\cal L}({\bf K}^{(n)})= {\cal L}(Z_1,\ldots, Z_n \vert \sum_{i=1}^niZ_i=n),\quad n=1,2,\ldots, \la{st1}
\end{equation}

\nin where $Z_1,Z_2,\ldots$ are independent integer valued random variables.
  The great importance  of  the  conditioning relation
\refm[st1] is based on the following two
  interrelated facts that hold for a variety of instances of RCS's.
\begin{itemize}

\item  The distribution  of $Z_i, \ i=1,2,\ldots$
is of one  of the following three types:

(i) Poisson($\frac{m_ix^i}{i!}, x>0$), (ii) Negative
binomial($m_i, x^i,\ x\in (0,1)$)

or (iii) Binomial($m_i,\frac{x^i}{1+x^i}, \ x>0$), where in all the cases $x$ is a free parameter and $m_i$ is the
number of components of size $i$.

\item Corresponding to the type of the distribution of $Z_i,$
the relationship between the two key sequences $\{p_n\}$ and $\{m_i\}$ has the form:

(i)  \be \sum_{n=0}^\infty \frac{p_nz^n}{n!}= \exp\Big(\sum_{i=1}^\infty\frac{m_iz^i}{i!}\Big), \la{st2}
\end{equation}

(ii) \be \sum_{n=0}^\infty p_nz^n=\prod_{i=1}^\infty (1-z^i)^{-m_i}, \la{st3}
\end{equation}

(iii) \be \sum_{n=0}^\infty p_nz^n=\prod_{i=1}^\infty (1+z^i)^{m_i}. \la{st4}
\end{equation}

\end{itemize}

\nin In accordance with the above, the following three basic classes of
  CS's are distinguished (\cite{abt}):

\nin  (i) assemblies, (ii) multisets and  (iii)
  selections.

\nin First, we immediately see  from \refm[st2] that assemblies are incorporated into our setting \refm[00] with
$a_n, c_n$
 having a clear combinatorial
context: $a_n=\frac{m_n}{n!}, \quad c_n=\frac{p_n}{n!}, \quad n=1,2,\ldots.$

\nin A quite different approach leading to the relationship \refm[st2] is widely known in  combinatorics (for
references see \cite{st}, Ch.5). In this field, \refm[st2] which is called the exponential formula, expresses the
general enumerative principle for  posets, that may be  regarded as disjoint unions of their connected components.
In particular, $S(z)=\sum_{n=1}^\infty\frac{m_nz^n}{n!}$ and $g(z)=\sum_{n=0}^\infty \frac{p_nz^n}{n!}$  are
called exponential generating functions for the number of connected components and for the total number of posets,
respectively. Note that in the graph theory, \refm[st2] is deduced from the generalized scheme of allocation (see
\cite{kol}, Ch.1, by V. Kolchin), the latter being  equivalent, in effect, to the aforementioned enumerative
principle.

 \nin   Multisets can  be also put into the framework of
\refm[00], by exponentiation of the generating function for
 the sequence $\{m_i\}$ (for references see Ch.2 of the
monograph \cite{bur}, by S.Burris). \nin We write \be \prod_{i=1}^\infty (1-z^i)^{-m_i}=\exp\Big(
\sum_{i=1}^\infty m_i \log(1-z^i)^{-1}\Big)= \exp\Big(\sum_{n=1}^\infty z^n\sum_{j,k:jk=n}\frac{m_j}{k}\Big)
\la{st5}
\end{equation}
\nin to get from \refm[st3],  $c_n=p_n, \quad a_n=\sum_{j,k:jk=n}\frac{m_j}{k}.$

\nin  Thus, the counting processes ${\bf{K}}^{(n)}$ for assemblies and multisets satisfy

\be {\cal L}({\bf K}^{(n)})=\mu_n, \la{st6}
\end{equation}

\nin where $\mu_n$ is the measure given by \refm[c2] with
$\gamma=1$ and the parametric function $a$ is  as indicated above.
\nin Though, in the case of multisets, $a_n$ lacks a combinatorial
meaning, it turns out that, under a certain condition, the
asymptotic behaviors of the two sequences $\{a_n\}$ and $\{m_j\}$
are similar.

\begin{prop} (\cite{bbur}, Lemma 5.1)

\nin If the sequence $\{m_j\}$ in \refm[st5] is such that \be \lim_{j\to \infty} \frac{m_j}{m_{j+1}}=h, \quad 0<h<1,
\la{st6*}
\end{equation}
\nin then $a_j\sim m_j, \quad n\to\infty.$
\end{prop}
\vskip .5cm \nin By virtue of  Remark 5 this means that our results on clustering are applicable for multisets
with $m_j\sim h^jj^{l-1} L(j), \ j\to \infty, \quad l,h>0.$
  Now notice that applying the exponentiation procedure in the case of selections we arrive at an alternating
sequence $\{a_n\}$ . This says that this case is beyond  the scope of the setting of
 the present paper.

\nin The asymptotic behavior of  counting processes
 was fully explored for the subclass of RCS's characterized by
 the following logarithmic condition:

 \be
\lim_{i\to \infty} i P(Z_i=1)=\lim_{i\to \infty} i EZ_i=\theta, \la{st7}
\end{equation}

\nin for some $\theta>0,$ where the random variables  $Z_i, \quad
i=1,2,\ldots $ are as in \refm[st1]. Such RCS's are called
logarithmic. A comprehensive exposition of the research for this
case is given in \cite{abt}.

\nin The  classical example of a logarithmic RCS is the seminal Ewens sampling formula (ESF) given by
$a_n=\frac{\theta}{n}, \ n\ge 1, \ \theta>0.$ It originated  in population genetics (1972) and was extensively
investigated by many authors in relation to a variety of models. In particular, it was proved that the normalized
ESF converges weakly to the Poisson - Dirichlet law (see \cite{ew}, \cite{hr} and references therein).
 The
counting process induced by ESF can be interpreted as a $theta-$ biased random permutation (\cite{abt}, Ch.3).
 The theory of the
limiting behavior of the counting process in the case $\theta=1$
(= random permutations) was shaped by V. L. Goncharov (1942), L.
A. Shepp  and S. P. Lloyd (1966) and A. M. Vershik and A. A.Shmidt
(1977) (for references see \cite{abt}, Ch.1 and \cite{kol}, Ch.4).

\nin On the other hand, integer partitions provide an example of a nonlogarithmic RCS. Partitions  can be formally
defined as a multiset with $m_i=1, \quad i\ge 1.$ Thus, \refm[st5] gives for this case
$$ a_n=\sum_{d\in D_n}\frac{1}{d}\le \sum_{j=1}^{\frac{n}{2}} \frac{1}{j},$$
\nin where $D_n$ is the set of all divisors of $n.$ Consequently, $$1\le a_n\le \log n, \quad n\to \infty,$$ \nin
which indicates that the case of partitions can be approximated by  the class of parametric functions ${\cal F}_l$
with $l=1$.

\nin In the next section we explain that $q$ -  colored linear
forests (treated as posets)  is a RCS  with $a\in {\cal F}_1.$ In
the conclusion, we make the following

\begin{rem} \label{p9}  The logarithmic condition \refm[st7] fails for the class
${\cal F}_l, \ l\neq 0,$
 of parametric functions $a$. On the other hand, the Lyapunov condition, and
 consequently the normal local limit  theorem, hold only when  $l>0.$
 This explains why   in the study of the clustering problem,
  the cases $l=0$, $l>0$ and $l<0$ should be distinguished,
  with basically different asymptotic tools being employed.
The third case that includes such RCS's , as forests of labelled
 (unlabelled) trees, was recently explored in \cite{bg}, by A.
Barbour and B. Granovsky.  Correspondingly, three very different
pictures of clustering were discovered. A specific feature of
clustering in the case $l>0$ considered in the present paper is
the existence of a threshold value for the size of the maximal
component (cluster). So, this appears to be  the only case (among
$a\in {\cal F}_l$)  in which the gelation phenomena is not seen.

\nin In the context of RCS's, the  aforementioned principle of
asymptotic independence of small groups (= components) has been
widely discussed for a long time, in connection with the
conditioning relation \refm[st1]. The independence  was proved for
the logarithmic RCS's (\cite{abt}, Ch. 4) and in  the case
$EZ_j=j^{l-1}L(j), \ l\le 0$  in \cite{bg}.
\end{rem}
 \nin \section{Application
3: Additive number systems (ANS).}

\nin ANS's provide a very general setting that encompasses multisets, as defined in the previous section.
 Following  \cite{bur} by S. Burris, an ANS ${\cal A}$ is a countable free
 commutative monoid $ A=\{v\}$ with a
given set $P$ of nondecomposable elements (= generators) and with
an additive norm  $\parallel\bullet\parallel,$ such that the set
 $$ \{v\in A:\parallel v
\parallel=n\} $$ \nin is finite for all $n\in {\cal N}.$ This
definition implies that each $v\in A$ is a sum (= union) of
elements of $P$. Denoting $c_n, \ p_n$ the  number of elements in
$A$ and $P$ correspondingly with norm $n,$ an enumerative argument
yields the following characteristic identity for ANS's: \be
\sum_{n\ge 0} c_n x^n=\prod_{n\ge 1}(1-x^n)^{-p_n},\quad 0\le x\le
\rho <1. \la{an1}
\end{equation}

\nin By the exponentiation of the RHS of \refm[an1], we get the
alternative version of the above identity: \be
g(x)=exp\Big(\sum_{m\ge 1}\frac{P(x^m)}{m}\Big), \quad 0\le x\le
\rho <1, \la{an2}
\end{equation}

\nin where $g$ and $P$ are the generating functions for the sequences $\{c_n\}$ and $\{p_n\}$ respectively: \be
g(x)=\sum_{n\ge 0} c_nx^n, \quad P(x)=\sum_{n\ge 0} p_n x^n, \quad 0\le x\le \rho <1. \la{an3}
\end{equation}

\nin Now \refm[an2] can be rewritten as \refm[00] with \be a_n=\sum_{jm=n}\frac{p_j}{m}. \la{an4}
\end{equation}

\nin As we already mentioned before, the sequence $\{ a_n\}$ defined by
\refm[an4] usually does not exhibit a  regular asymptotic behaviour, i.e. does not belong to the class  ${\cal
F}(l), \  l> 0$.

\nin On the other hand, ANS's with $p_n$ satisfying the condition
of Proposition 3 are a nice  exception to the above phenomenon. An
example of such a structure is the set of $q$-colored linear
forests, treated as posets, in which case $p_n=q^n$ (see
\cite{bur}, p.24).
 It is also important to note that
 the radius of convergence of the generating series in \refm[an3]
cannot be greater than $1.$

\nin  We wish now to  demonstrate  that the known asymptotic
result on ANS's that facilitated the development of Compton's
density theory  is a particular case of our asymptotic formula for
$c_n$. In 1992 it was proven by J. Knopfmacher, A. Knopfmacher and
R. Warlimont (see for references \cite{bur}, Theorem 5.17, p.94)
that if in \refm[an1], $$p_n=hq^n+O(q_1^n), \quad h>0, \quad
q>1,\quad 0<q_1<q,$$

\nin then \be c_n\sim h_1 q^n\frac{e^{2\sqrt{hn}}}{n^{\frac{3}{4}}}, \quad h_1>0,\quad n\to
\infty,\la{an5}\end{equation}

\nin where $h_1>0$ is a  constant which was not specified.
 This result was obtained with
the help of complex analysis.

\nin By virtue of Proposition 3, we see that $a_n\sim p_n, \ n\to
\infty,$ which together with Remark 5, permits to apply our
asymptotic formula \refm[836] with $l=1$ and $ L(n)\equiv h.$ So,
in the case considered, $L^*(n)\equiv h^{-1}$ and we have in
\refm[*], $L_1(n)=h^{-1/2}.$ Consequently, by the asymptotic
formulae in Section 3, we have $\un{\sigma}_n^{(n)}\sim
n^{-1/2}h^{1/2},$  $\un{S}_n^{(n)}(e^{-\un{\sigma}_n^{(n)}})\sim
h^{1/2}n^{1/2}-h/2+O(\un{\sigma}_n^{(n)}),$  and
$(\un{B}_n^n)^2\sim 2h^{-1/2}n^{3/2},$ as $n\to \infty.$
  Substituting  this in \refm[836], recovers \refm[an5],
while specifying $h_1=(2\sqrt{\pi})^{-1}h^{1/4}e^{-h/2}.$

\nin The central problem in the theory of ANS's is the study of the asymptotic density $\delta(B)$ of a subset
$B$ of a monoid $A$: \be \delta(B)=\lim_{n\to \infty} \frac{b_n}{c_n}, \la{an6}
\end{equation}
\nin where $b_n$ is the number of elements of $B$ with norm $n$.

 \nin It follows from \refm[an4] that
the quantities $ \underline{d}^{(r)}_n, \ \bar{d}^{(r)}_n$ in the clustering problem  considered in
the present paper can be regarded as the densities of the subsets, say $B_1,B_2\subseteq A$,
      such
 that the
  maximal (minimal) norm of  generators of  elements of $B_1 \ ( B_2)$
  satisfies a certain condition.

\nin The fundamental result in this field is Compton's density theorem (1989)
 (see \cite{bur}, Ch.4, 5) that
establishes sufficient conditions  for existence of an asymptotic
density of all partition sets of an ANS ${\cal A}.$

\nin Coming back to the example of $q$-colored linear forests (as posets), our Theorem 2, applied with $l=1$ and
$L(n)\equiv h$, gives the asymptotic density of the aforementioned sets $B_1,B_2$.

 \vskip 2cm
 \nin {\bf
Acknowledgement}

 The research of the second author was supported by E. and
 J. Bishop Research Fund.

The referee's criticism and remarks contributed to the improvement
of the exposition of

the paper.

\end{document}